\documentstyle[preprint,aps,prc]{revtex}
\begin{document}
\renewcommand{\baselinestretch}{2}
\tightenlines
 
\title{An accurate and efficient algorithm for the computation
 of the characteristic polynomial of a general square matrix.}
\author{S. Rombouts \thanks{Postdoctoral Fellow of the Fund for 
                            Scientific Research - Flanders (Belgium)} 
        and K. Heyde}
\address{Institute for Theoretical Physics,
         Department of Subatomic and Radiation Physics.
         \\
         Proeftuinstraat 86, B-9000 Gent,  Belgium
         \\
         E-mail: Stefan.Rombouts@rug.ac.be, Kris.Heyde@rug.ac.be
         }
\date{September 25, 1997}
\maketitle

\begin{abstract}
\noindent  
An algorithm is presented for the efficient and accurate 
computation of the coefficients of the characteristic polynomial of 
a general square matrix. 
The algorithm is especially suited for the evaluation of canonical traces
in determinant quantum Monte-Carlo methods.
\end{abstract}

\section{introduction}
The characteristic polynomial $P_U \left( x \right)$ 
of a general $N$-by-$N$ matrix $U$ is given by 
\begin{equation}
\label{charpol1}
P_U \left( x \right) = \det \left(x I - U \right),
\end{equation}
($I$ is the unit matrix).
Though the characteristic polynomial of a matrix is a basic concept 
in linear algebra, its numerical computation is scarcely documented 
in literature. 
There are however a number of applications for which an accurate and efficient 
algorithm for the calculation of the coefficients of $P_U \left( x \right)$
would be useful, 
e.g. for the study of random matrices \cite{haake}.
It would also be useful
for determinant quantum Monte-Carlo (QMC) methods \cite{gfqmc,hirsch,smmc}:
the application of determinant QMC method in the canonical ensemble,
especially the shell-model Monte-Carlo method \cite{smmc}, 
requires the evaluation of the coefficients of the polynomial
\begin{equation}
\label{charpol2}
  \bar{P}_U \left( x \right) = \det \left( I + x U \right).
\end{equation}
This polynomial is closely related to $P_U \left( x \right)$:
\begin{equation}
\label{charpolrel}
  \bar{P}_U \left( x \right) = 
     \left(-x\right)^N P_U \left(-1/x \right).
\end{equation}
In the determinant QMC method,
the $A$-particle trace of a one-body evolution matrix $U$ is given by
the coefficient of $x^A$ in $\bar{P}_U \left( x \right)$.
This coefficient is equal to $\left(-1\right)^A$ times the coefficient 
of $x^{\left(N-A \right)}$ in $P_U \left( x \right)$. 
It is in the light of these canonical QMC methods
that we developed an algorithm that is presented in the next sections.
Accuracy is important here because for calculations at low temperature,
QMC methods tend to become very sensitive 
to numerical instabilities in the evaluation of the matrix elements 
and the coefficients of $\bar{P}_U \left( x \right)$ \cite{lind}. 
Speed is important here because for one Monte-Carlo run,
several thousand matrices have to be evaluated.

The algorithms that can be found in literature,
as e.g. the Faddeev-Leverrier method \cite{faddeev}, 
are far from optimal on the point of view of numerical computation.
Let $a_A$ denote the coefficient of $x^A$ in $P_U \left( x \right)$
and let $t_k=Tr \left(U^k\right)$.
The coefficients $t_k$ correspond to the power sums of the roots 
of $P_U \left( x \right)$.
Their computation is straightforward.
A relation between the $a_A$ and the $t_k$ can be obtained by 
equating the coeffients of the powers of $x$ in the series expansion
of both sides of
\begin{equation}
  \det \left(x I - U \right) 
    = x^N \exp \left[ \mbox{Tr} \ln \left(I-\frac{1}{x} U\right) \right].
\end{equation}
This leads to a relation already established by Newton \cite{newton},
that formally can be written as 
\begin{equation}
  a_A = \frac{1}{A!} \det \left(
  \begin{array} {cccccc}
	t_1 &   1 &   0 &   0 & \cdots &   0 \\
	t_2 & t_1 &   2 &   0 & \cdots &   0 \\
	t_3 & t_2 & t_1 &   3 & \cdots &   0 \\
	\cdots & \cdots & \cdots & \cdots & \cdots & \cdots \\		       
	t_A & t_{A-1} & t_{A-2} & t_{A-3} & \cdots &   t_1 \\
  \end{array} \right).
\end{equation}
The Faddeev-Leverrier method \cite{faddeev} is a recurrence relation
that implicitly generates the same relations.
Though mathematically elegant, this formula is unpractical: 
it is inaccurate because it is very sensitive to roundoff errors, 
especially if the eigenvalues of $U$ differ 
by several orders of magnitude, which is common in determinant QMC methods. 
It is also inefficient because it requires $A/2$ matrix multiplications. 
Therefore the method is only useful for small $A$.
Another method, suggested by Ormand et al. \cite{ormand},
amounts to the following expression:
\begin{equation}
  a_A = \frac{1}{N} \sum_{m=0}^{N-1} e^{- i \frac{2 \pi}{N} m A} 
	\det \left( e^{i \frac{2 \pi}{N} m } I - U \right).
\end{equation}
In order to evaluate this expression efficiently,
it is suggested to diagonalize $U$ first.
However, if $U$ is diagonalized, 
the coefficients of $P_U \left( x \right)$ can be evaluated more easily 
by explicit construction of the polynomial
\begin{equation}
\label{detprod}
\det \left( x I - U \right) 
  = \prod_{i=1}^N \left( x - \epsilon_i \right) .
\end{equation}
If this polynomial in $x$ is constructed 
from the smallest up to the largest eigenvalue, 
$a_A$ can be computed in an easy and stable way. 
The main computational effort is in the diagonalization of the matrix $U$.
The polynomial can be constructed even more efficiently 
without diagonalization, as is explained in the next section.

In section II we present an algorithm for calculating the coefficients 
of the characteristic polynomial of a general square matrix. 
In section III we present the results of the numerical tests 
of the speed and accuracy of the algorithm.
 
\section{Algorithm for the calculation of the characteristic polynomial 
         of a general square matrix}

The basic idea of the algorithm is to consider $I+x U$ 
as a matrix of polynomials in $x$. 
We then calculate the polynomial $\bar{P}_U \left( x \right)$ by
evaluating the determinant in equation \ref{charpol2} 
using Gaussian elimination, 
with polynomials instead of scalars as matrix elements. 
As mentioned above, the coefficients of $\bar{P}_U \left( x \right)$ 
are closely related to the coefficients $a_A$.
Because the multiplication of two polynomials of degree $A$ 
requires about $2A^2$ flops and calculation of a determinant 
about $N^3/3$ polynomial multiplications, 
the calculation would require a number of the order of $N^5$ flops, 
which is too much for an efficient implementation. 
This number can be drastically reduced if $U$ is transformed 
to an upper-Hessenberg form by a similarity transformation 
(a Householder reduction to Hessenberg form requires 
approximatly $\frac{10}{3}N^3$ flops \cite{matrix}). 
This leaves the coefficients of $\bar{P}_U \left( x \right)$ unchanged.
In order to calculate the determinant we transform $I + x U$ 
to upper diagonal form by Gaussian elimination, 
requiring now only $N^2$ polynomial multiplications. 
The Gaussian elemination is performed from the right bottom corner 
of the matrix up to the top left corner 
because in determinant QMC methods, 
the right bottom corner often contains the smallest elements,
so that the summations involved are performed from small to large terms, 
which is less sensitive to roundoff errors 
than the summation the other way round.
We start with $T^N$=$I+x U$. 
Now we bring column after column in upper triangular form.
Suppose that $T^j$ has column $j$ to $N$ already in upper triangular form, i.e.
\begin{equation}
\label{uptriang}
T^j_{i \, k} = 0,
\end{equation}
for $i>k$ and $k \geq j$.
Now we calculate 
\begin{equation}
\label{timesg}
T^{j-1} = T^{j} G^{j},
\end{equation}
where
\begin{equation}
\label{gis1}
G^{j}_{i \, k} = \delta_{i,k},
\end{equation}
except for
\begin{equation}
\begin{array} {lcl}
\label{gis2}
G^j_{j-1 \, j-1} & = &  T^j_{j \, j} , \\
G^j_{j \, j-1}   & = & -T^j_{j \, j-1}.
\end{array}
\end{equation}
In the end we obtain the upper triangular matrix 
$T^1= T^N \, G^N \, G^{N-1} \, \cdots \, G^{2}$ so that
\begin{eqnarray}
\label{tprod}
\bar{P}_U \left( x \right) & = & \det \left( I + x U \right) \\
 & = & \det \left( T^N \right) \\
 & = & \frac{\det \left( T^1 \right)}
       {\det \left(  G^N \, G^{N-1} \, \cdots \, G^{2} \right)} \\
 & = & \frac{\prod_{i \, =1}^N T^1_{i \, i}}{\prod_{i \,=2}^N T^i_{i \, i}} \\
 & = & T^1_{1 \, 1}.
\end{eqnarray}
because $T^1_{i \, i} = T^i_{i \, i}$.
The operations can be ordered to minimize memory use.
This leads to the following algorithm 
($t_{k \,i}$ corresponds with the coefficient of $x^k$ in $T^j_{i \, j}$) :
\begin{equation}
\label{algo1}
\nonumber
\begin{array}{l} 
{\bf algorithm \; for \; calculating \; the \; coefficients \; 
 of \; the \; characteristic} \\
{\bf  polynomial \; of \; a \; N-by-N matrix \; U} \\
 reduce \; U \; to \; upper Hessenberg \; form	  \\
 \mbox{DO} \ j=N,1,-1                             \\ 
 \; \; \; \; \; \;  
 \mbox{DO} \   i=1,j                              \\ 
 \; \; \; \; \; \; \; \; \; \; \; \;    
 \mbox{DO} \   k=N-j,1,-1                         \\ 
 \; \; \; \; \; \; \; \; \; \; \; \; \; \; \; \; \; \; 
 t_{k+1 \, i}=U_{i \, j} \, t_{k \, j+1}-U_{j+1 \, j} \, t_{k \, i}  \\ 
 \; \; \; \; \; \; \; \; \; \; \; \;
 \mbox{ENDDO}                                     \\ 
 \; \; \; \; \; \; \; \; \; \; \; \;
 t_{1 \, i}=U_{i \, j}                            \\ 
 \; \; \;  \; \; \; 
 \mbox{ENDDO}                                     \\ 
 \; \; \; \; \; \;
 \mbox{DO} \  k=1,N-j                             \\
 \; \; \; \; \; \; \; \; \; \; \; \;
 t_{k \, j}=t_{k \, j}+t_{k \, j+1}               \\ 
 \; \; \; \; \; \;
 \mbox{ENDDO}                                     \\ 
 \mbox{ENDDO} 
\end{array}
\end{equation}
In the end $t_{k \, 1}$ is the coefficient of $x^k$ in $\bar{P}_U \left( x \right)$.
Then $a_A$ is given by
\begin{equation}
a_A = \left(-1\right)^{N-A} t_{N-A \; 1}.
\end{equation}
This algorithm cannot break down and requires $N^3 /2+N^2-N/2$ flops. 
If one needs only the coefficient of $x^A$ in $\bar{P}_U \left( x \right)$, 
e.g. for the calculation of an $A$-particle trace in determinant QMC methods, 
the number of flops can be reduced further by calculating 
the polynomials only up to degree $A$. 
This is done by
restricting the loop over k to values smaller than or equal to $A$.
The sixth line in \ref{algo1} then becomes   
\begin{equation}
\label{krest}
  \mbox{DO} \  \ k= \mbox{MAX} \left(N-j,A \right) ,1,-1 .
\end{equation}
Together with the Householder reduction to the upper Hessenberg form 
this makes less than $4 N^3$ flops. 
Diagonalization of the matrix $U$ with the QR algorithm (suggested in 
\cite{matrix} as the obvious method for the diagonalization of general 
square matrices), would require about $10 N^3$ flops. 

\section{Numerical tests}
We have tested our algorithm numerically on its speed and accuracy.
All the tests were done in Fortran77 (DEC Fortran V3.8) on 
a Digital Alphastation 255/300MHz workstation running Digital Unix 3.2D. 
For the reduction to hessenberg form and the diagonalization 
optimized Lapack routines were used \cite{lapack}. 
For the part of the algorithm listed in the previous section 
only the standard optimizations of the Fortran compiler were used.  

The speed was tested by calculating, for several matrix sizes, 
all the coefficients of the characteristic polynomial of 100 matrices 
with random elements (uniformely distributed 
between $-\frac{1}{\sqrt{N}}$ and $\frac{1}{\sqrt{N}}$). 
This was done with our algorithm and with complete diagonalization. 
The speed was measured by counting the number of cycles 
executed by the procedures of the algorithms 
(fewer cycles means faster calculation) using the 'prof -pixie' command.
Table I lists the results. 
It is clear that our algorithm is much faster than complete diagonalization: 
from a factor 4.5 for small matrices to a factor 1.8 for large matrices. 
The decrease of this factor for large matrices can be understood 
by the fact that the routines for the reduction to hessenberg form 
and diagonalization are strongly optimized while the routine 
for the algoritm of section III is not, 
and that these optimizations become more and more efficient 
with larger matrix sizes.
Applying similar optimizations to our algorithm would
significantly increase the ratio.  
Typical matrix sizes in determinant QMC methods range from 20-by-20 
to 100-by-100. 
So there our unoptimized algorithm is a factor 2 to 4 faster
than complete diagonalization.  

For testing the accuracy, we calculated 200000 random samples 
with a Determinant Quantum Monte-Carlo method 
for the 4x4 Hubbard model with 8 up and 8 down electrons,
with U=4 and $\beta$=6, following the method of reference \cite{white}, 
but taking the canonical trace instead of the grand-canonical one. 
For each sample, the canonical trace is given by the square of the 
coefficient of $x^8$ in the characteristic polynomial 
of a 16-by-16 matrix.
The canonical trace was calculated in double precision and in single precision 
using our algorithm and complete diagonalization. 
As a measure for the accuracy we used the average absolute value 
of the difference between the single- and double-precision result 
divided by the double-precision result. 
For our algorithm we found a value of 0.00186 $\pm$ 0.00005 
and for the complete diagonalization we found 0.00607 $\pm$ 0.00010 
(error limits at 95\% confidence level), 
indicating that our algorithm is more accurate. 
This could be expected since it requires less operations on the data. 
Furthermore complete diagonalization was much more sensitive 
to overflow errors than our algorithm. 
At values of $\beta > 6$ complete diagonalization 
(in single precision) was not usable anymore.  

\section{Conclusion}
We have presented a stable and efficient 
algorithm for the calculation of the coefficients 
of the characteristic polynomial of a general square matrix. 
This algorithm is especially useful for
Determinant quantum Monte-Carlo calculations in the canonical ensemble
because it is faster (a factor 2 to 4) and more accurate 
than the algorithms that can be found in the literature. 

\section*{Acknowledgements}
This work was supported by the 
Fund for Scientific Research - Flanders (Belgium).

\renewcommand{\baselinestretch}{1}
\begin{table}
\begin{tabular}{|r ||r |r |r |}
{\em matrix dimension} & {\em our algorithm} 
                       & {\em complete diagonalization} & ratio \\
\hline
     4 &          451400  &        1983818   &   4.39 \\ 
     6 &         1009400  &        4413843   &   4.37 \\
     8 &         1760300  &        8062663   &   4.58 \\
    10 &         2870100  &       12511637   &   4.36 \\
    15 &         7261300  &       29676436   &   4.09 \\
    20 &        14224100  &       55696656   &   3.93 \\
    25 &        25524300  &       93696774   &   3.67 \\
    30 &        41735900  &      144177197   &   3.45 \\
    35 &        63852100  &      209670202   &   3.28 \\
    40 &        90395400  &      290658105   &   3.22 \\
    45 &       126513800  &      388488344   &   3.07 \\
    50 &       171095900  &      512056714   &   2.99 \\
    60 &       284484900  &      794032492   &   2.79 \\
    70 &       447113900  &     1163945220   &   2.60 \\
    80 &       652709400  &     1630550332   &   2.50 \\
    90 &       926006900  &     2207209655   &   2.38 \\
   100 &      1251268900  &     2923248380   &   2.34 \\
   150 &      4268580600  &     8925077120   &   2.09 \\
   200 &     10018384500  &    20050929483   &   2.00 \\
   300 &     32993383700  &    63384810388   &   1.92 \\
   400 &     77249914100  &   145321243773   &   1.88 \\
   500 &    149825926400  &   278218705522   &   1.86 \\
   600 &    257427888100  &   474763616745   &   1.84 \\
   700 &    407433443000  &   745631287828   &   1.83 \\
   800 &    607094132500  &  1104878051129   &   1.82 \\
   900 &    863225666500  &  1564619645628   &   1.81 \\
\hline
\end{tabular}
\caption{Comparision of the number of cycles needed for the calculation 
         of the coefficients of the characteristic polynomial of 100 
         matrices with random elements, for several matrix dimensions.}
\end{table}

\end{document}